# When is $\sin x + \cos x + \tan x + \cot x + \sec x + \csc x$ an integer ?

*Konstantine Zelator*
*Department Of Mathematics*
*College Of Arts And Sciences*
*Mail Stop 942*
*University Of Toledo*
*Toledo,OH 43606-3390*
*U.S.A.*

## 1. Introduction

Consider the function,
*F(x)=sinx+cosx+secx+cscx+tanx+cotx,*

with domain all real numbers which are not of the form $K\pi$ or $K\pi + \frac{\pi}{2}$;

K an integer. This is a periodic function of period $2\pi$.
The motivating problem behind this work, is found in Section 3.
The original book source for that motivating problem, was an obscure trigonometry book (in Greek) published in Athens, Greece, in 1973 (see [1]). Currently, this problem can also be found in [2].
 In Section 2 of this paper we list six very widely and extremely well known results, mostly trigonometric identities from trigonometry – no further commentary needed. We use these results in Sections 3 and 5. In Section 4, we state four results from classical algebra, and we offer short proofs for two of them.
Section 5 is the central focus of this article. We show that each integer n which satisfies n ≥ 7 or n ≤ -2 lies in the range of the above function F. For any such integer n, we find and explicitly state the pre-images of n in the domain of F; in other words, we find the solutions to the equation F(x) =n. By contrast, as we will demonstrate, the latter equation has no solution for n = -1, 0, 1, 2, 3, 4, 5, or 6.

## 2. Six results from Trigonometry

**Result 1**: For any values of the radian (or degree) angle measure x,
$$\sin^2 x + \cos^2 x = 1$$

**Result 2**: For any values of the radian (or degree) angle measures A and B,
$$\cos A + \cos B = 2\cos\left(\frac{A+B}{2}\right)\cos\left(\frac{A-B}{2}\right)$$

**Result 3**: For any value of the radian (or degree) angle measure of t,
$$\cos(-t) = \cos t$$

**Result 4**: For any value of the radian (or degree) angle measure of x,
$$\sin x = \cos\left(\frac{\pi}{2} - x\right)$$

**Result 5**: For any value of the radian (or degree) angle measure of x,
$$S = \sin x + \cos x = \sqrt{2}\cdot\cos\left(x - \frac{\pi}{4}\right); \text{ and} -\sqrt{2} \leq S \leq \sqrt{2}.$$

Observe that
$$S = \sin x + \cos x = \cos(\pi/2 - x) + \cos x = 2\cos\left[\frac{(\pi/2 - x) + x}{2}\right]\cos\left[\frac{x - (\pi/2 - x)}{2}\right];$$

$$S = 2\cos\pi/4 \cos(x - \pi/4) \overset{\text{by Result 4}}{=} 2 \cdot \frac{\sqrt{2}}{2} \cos(x - \pi/4) \overset{\text{by Result 2}}{=} \sqrt{2}\cos(x - \pi/4)$$

**Result 6**: Consider the fundamental trigonometric equation $\cos x = b$.

(i) If $|b| > 1$; that is if b>1 or b<-1; the equation has no solutions.
(ii) If b=-1, all the solutions are given by $x = (2K+1)\pi$, where K can be any integer.
(iii) If b=1, all the solutions are given by $x = 2K\pi$; K can be any integer.
(iv) If -1<b<1; let $\theta_b$ be the unique radian angle measure such that $0 < \theta_b < \pi$ and $\cos\theta_b = b$. Then, all the solutions to the above equation are given by $x = 2K\pi \pm \theta_b$, where K can be any integer.

### 3. The motivating problem and its solution

In 1948, the following problem was given in the entrance exams for the National Technical University of Athens, Greece. In the years following World War II, those entrance exams had become notoriously difficult in mathematics (algebra, geometry, and trigonometry). Now the problem:
Find all the solutions of the trigonometric equation,

$$\sin|x| + \cos|x| + \tan|x| + \frac{1}{\sin|x|} + \frac{1}{\cos|x|} + \frac{1}{\tan|x|} + 3 = 0 \quad (1)$$

In effect, equation (1) is equivalent to,
$$F(|x|) = -3.$$

Obviously, any solution x to (1) requires $|x| \neq K\pi, \ K\pi + \pi/2$; K an integer.

The technique we use to solve (1), involves the sum $S = \sin|x| + \cos|x|$.

From $S^2 = \sin^2|x| + \cos^2|x| + 2\sin|x|\cos|x|$ and Result 1, we easily obtain,

$$\sin|x|\cos|x| = \frac{S^2 - 1}{2} \quad (2)$$

Note that because of the restriction $|x| \neq K\pi, K\pi + \pi/2$; we must have $S \neq 1, -1$.

If we multiply equation (1) by $\sin|x|\cos|x|$ and use (2) we obtain,

$$S(S^2 - 1) + 2S + 2 + 3(S^2 - 1) = 0;$$
$$S^3 + 3S^2 + S - 1 = 0 \quad (3)$$

The number -1 is a root to equation (3). Using synthetic division (or long division) we obtain the factorization,

$$(S + 1)(S^2 + 2S - 1) = 0 \quad (4)$$

Thus, the other two roots of equation (3) are the two zeros of the trinomial $S^2 + 2S - 1$; which are the real umbers $(-1 + \sqrt{2})$ and $-(1 + \sqrt{2})$ (by use of quadratic formula). Of these three roots; -1 is ruled out, since $S \neq -1, 1$, as we have explained above. Also, $-\left(1+\sqrt{2}\right)$ does not fall in the interval $\left[-\sqrt{2},\sqrt{2}\right]$ as required by Result 5 ($S$ can only assume values between $-\sqrt{2}$ and $\sqrt{2}$). Therefore,

$S = -1+\sqrt{2}$; And by Result 5, $\cos(|x|-\pi/4) = -\dfrac{1}{\sqrt{2}}+1 > 0$.

In the language of Result 6, $b = 1-\dfrac{1}{\sqrt{2}}$; $-1 < b < 1$. Let $\varphi$ be the unique radian angle measure such that $\cos\varphi = 1-\pi/\sqrt{2}$ and $0 < \varphi < \pi$; in fact $0 < \varphi < \pi/2$, since $1-1/\sqrt{2} > 0$. According to Result 6(*iv*), we arrive at $|x| = 2K\pi \pm \varphi$; $\boxed{x = \pm(2K\pi \pm \varphi)}$; K can be any integer; with the choice of signs being arbitrary (so, all together, there are four sign combinations).

Note that $\varphi \approx 1.27354433$ radians ($\approx 72.9687153$ degrees).

The above described problem can be found in reference [2]. However, the original source of the problem is [1].

**4. Four results from classical algebra.**

The first of these results is the factorization of a polynomial of degree n with real coefficients into n linear factors involving its n roots or zeros.

The said result can be found in college algebra and precalculus texts. Here, we are only interested in the case n=2.

---

**Result 7**. If $f(x) = ax^2 + bx + c$ is a quadratic function (i.e. $a \neq 0$) with real coefficients then,
$f(x) = a(x-p_1)(x-p_2)$, where $p_1$ and $p_2$ are the two zeros of f.

---

Next, if we complete the square we find that,

$f(x) = a\left(x+\dfrac{b}{2a}\right)^2 + \dfrac{4ac-b^2}{4a}$; which in turn, establishes the following.

---

**Result 8**: Let $f(x) = ax^2 + bx + c$; $a \neq 0$, b, c, being real numbers.
(*i*)   If $b^2 - 4ac < 0$ and $a > 0$, then $f(x) > 0$, for all real values of x.
(*ii*)  If $b^2 - 4ac < 0$ and $a < 0$, then $f(x) < 0$, for all real values of x.
(*iii*) If $b^2 - 4ac = 0$ and $a > 0$, then $f(x) > 0$ for any real value of $x \neq -b/2a$.
        And $f(-b/2a) = 0$.
(*iv*)  If $b^2 - 4ac = 0$ and $a < 0$, then $f(x) < 0$, for all real $x \neq -b/2a$.
        And $f(-b/2a) = 0$.

---

Results 9 and 10 below are typically not found in college algebra texts; certainly not in the explicit form presented below. For this reason, we offer a short proof for each of them.

**Result 9**: Let $\kappa$ and $\lambda$ be two real numbers, with $\kappa < \lambda$. Consider the trinomial $f(x) = ax^2 + bx + c$, with real coefficients and $a \neq 0$. Then f will have two distinct real roots both lying in the open interval $(\kappa, \lambda)$ if, and only if, all three conditions listed below are satisfied:
*Condition 1*: $b^2 - 4ac > 0$.
*Condition 2*: $af(\kappa) > 0$ and $af(\lambda) > 0$.
*Condition 3*: $\kappa < -\dfrac{b}{2a} < \lambda$.

**Proof.** We leave it to the reader to prove that the three conditions are necessary. We prove that they are sufficient. By Condition 1, the trinomial f must have two distinct real roots $p_1$ and $p_2$; say $p_1 < p_2$. By Condition 2 and Result 7, $a^2(\kappa - p_1)(\kappa - p_2) > 0$ and $a^2(\lambda - p_1)(\lambda - p_2) > 0$; and since $a^2 > 0$; it follows that $(\kappa - p_1)(\kappa - p_2) > 0$ and $(\lambda - p_1)(\lambda - p_2) > 0$. In each case, either both linear factors are positive, or both negative.

Thus, both $\kappa$ and $\lambda$ must lie outside the closed interval $[p_1, p_2]$. The midpoint of this interval is $\dfrac{p_1 + p_2}{2} = \dfrac{1}{2} \cdot \left(-\dfrac{b}{a}\right) = -\dfrac{b}{2a}$, since $p_1 + p_2 = -\dfrac{b}{a}$ (one way to establish this is straight from the quadratic formula). By Condition 3, then, the midpoint of the interval $[p_1, p_2]$ falls inside the open interval $(\kappa, \lambda)$. And since $\kappa < \lambda$ and both $\kappa$ and $\lambda$ lie outside $[p_1, p_2]$; we infer that $\kappa < p_1$ and $p_2 < \lambda$; $\kappa < p_1 < p_2 < \lambda$. □

**Result 10**: Let $\kappa$ and $\lambda$ be real numbers, $\kappa < \lambda$; and $f(x) = ax^2 + bx + c$, with a, b, c being real numbers; $a \neq 0$. Then, a necessary and sufficient condition for the trinomial f to have one real root inside the open interval $(\kappa, \lambda)$; while the other real root lying outside the closed interval $[\kappa, \lambda]$ (that is, in $(-\infty, \kappa) \cup (\lambda, +\infty)$); is $f(\kappa)f(\lambda) < 0$.

**Proof**: Again, as in the previous proof, we prove the sufficiency of the condition; the necessity is left to the reader. The condition $f(\kappa)f(\lambda) < 0$ is equivalent to,
($f(\kappa) < 0$ and $f(\lambda) > 0$) or ($f(\kappa) > 0$ and $f(\lambda) < 0$).
In either case, it follows from Result 8 that $b^2 - 4ac > 0$; that is, the trinomial f has two distinct real roots $p_1$ and $p_2$; say $p_1 < p_2$. Then,
$f(\kappa)f(\lambda) < 0 \Leftrightarrow a^2 \cdot f(\kappa)f(\lambda) < 0 \Leftrightarrow$ either ($af(\kappa) < 0$ and $af(\lambda) > 0$) or ($af(\kappa) > 0$ and $af(\lambda) < 0$).

In the first case, by applying Result 7, it follows that $a^2 \cdot (\kappa - p_1)(\kappa - p_2) < 0$ and $a^2 \cdot (\lambda - p_1)(\lambda - p_2) > 0$; which implies $p_1 < \kappa < p_2$ and $\lambda \in (-\infty, p_1) \cup (p_2, +\infty)$. In the second case, it is the other way around: $p_1 < \lambda < p_2$ and $\kappa \in (-\infty, p_1) \cup (p_2, +\infty)$. Thus, in the first case, $P_2 \in (\kappa, \lambda)$ and $P_1 \in (-\infty, \kappa)$; while in the second case, $P_1 \in (\kappa, \lambda)$ and $P_2 \in (\lambda, +\infty)$.  □

## 5. The equation F(x)=n, n an integer.

In this section, we find all the solutions to the equation,
$$F(x) = n;$$
$$\sin x + \cos x + \sec x + \csc x + \tan x + \cot x = n \qquad (5),$$
where n is an integer.

For convenience, we put $m = -n$. Also recall that $x \neq K\pi + \pi/2, \ K\pi$; K an integer. Putting all this together,

$$\begin{cases} m = -n \\ S = \sin x + \cos x; \\ \text{and so}, -\sqrt{2} \leq S \leq \sqrt{2} \\ \qquad \text{(by Result 5)} \\ \text{Also}, x \neq K\pi, K\pi + \pi/2; K \in Z \\ \text{Equivalently}, S \neq -1, 1 \end{cases} \qquad (6)$$

Also, we express the product sinxcosx in terms of S (by using Result 1):
$$\sin x \cos x = \frac{S^2 - 1}{2} \qquad (7)$$

By applying the standard identities
$\sec x = \frac{1}{\cos x}, \csc x = \frac{1}{\sin x}, \tan x = \frac{\sin x}{\cos x}, \cot x = \frac{\cos x}{\sin x}$; multiplying equation (5) with the product sinxcosx; and then using (6) and (7), and after some algebra we obtain the equivalent equation,
$$S^3 + mS^2 + S + (2 - m) = 0 \qquad (8)$$
The number -1 is a root of equation (8).
We perform synthetic division:

```
-1 |  1    m      1      (2-m)
              (m-2)
      1   (m-1)  (2-m)    0
```

Thus, equation (8) is equivalent to,
$$(S+1)[S^2 + (m-1)S + (2-m)] = 0 \qquad (9)$$
According to (6), $S \neq -1$. Therefore, by (9),
S must equal one of the two zeros of the trinomial $f(t) = t^2 + (m-1)t + (2-m)$.
Next, note that neither $\sqrt{2}$ or $-\sqrt{2}$ can be a zero of f. This inference follows easily from the fact that $\sqrt{2}$ is an irrational number, while m is an integer.
But then, according to (6), the sum S can be equal to one of the two zeros of the trinomial f; only if that zero falls within the open interval $(-\sqrt{2}, \sqrt{2})$. There are three possibilities:

Possibility 1: The trinomial f(t) has two distinct real roots, both of which fall in
$(-\sqrt{2}, \sqrt{2})$.

Possibility 2: The trinomial f(t) has two distinct real roots, one of which falls in
$(-\sqrt{2}, \sqrt{2})$; while the other lies outside $[-\sqrt{2}, \sqrt{2}]$; that is in
$(-\infty, -\sqrt{2}) \cup (\sqrt{2}, +\infty)$.

Possibility 3: The trinomial f(t) has a double real root in $(-\sqrt{2}, \sqrt{2})$.

Let us examine the third possibility first. It would require that the discriminant of the trinomial f(t) be zero:

$(m-1)^2 - 4(2-m) = 0 \Leftrightarrow m^2 + 2m - 7 = 0$, which is impossible is m is an integer and the two roots of the last quadratic equation are the irrational numbers $(-1-2\sqrt{2})$ and $(-1-2\sqrt{2})$.

Thus, Possibility 3 is ruled out. Next, consider Possibility 1. We apply Result 9 with $\kappa = -\sqrt{2}$ and $\lambda = \sqrt{2}$. The coefficients of the trinomial f(t) are $a = 1, b = m-1, c = 2-m$. Accordingly, the following four inequalities must simultaneously hold true:

$$\begin{cases} (m-1)^2 - 4(2-m) > 0, \\ \text{and} \quad 2 + (m-1)\sqrt{2} + 2 - m > 0; \\ \text{and} \quad 2 - (m-1)\sqrt{2} + 2 - m > 0; \\ \text{and} \quad -\sqrt{2} < -\frac{(m-1)}{2} < \sqrt{2} \end{cases}$$

$$\Leftrightarrow \begin{cases} m^2 + 2m - 7 > 0; \text{ and} \\ m > -\left(\frac{4-\sqrt{2}}{\sqrt{2}-1}\right) \approx -6.242640693; \text{ and} \\ m < \frac{4+\sqrt{2}}{1+\sqrt{2}} \approx 2.242640687; \text{ and} \\ 1-2\sqrt{2} \approx -1.828427125 < m < 1+2\sqrt{2} \approx 3.828427125 \end{cases} \quad (10)$$

Since m is an integer; by inspection we see that the last three inequalities in (10) are simultaneously satisfied exactly when m = -1, 0, 1, or 2. Of these four integers, only m = 2 satisfies the first inequality in (10). For m = 2; n = −2 (from (6)). The trinomial f(t) becomes, $f(t) = t^2 + t = t(t+1)$, which has zeros the numbers 0 and -1. But $S \neq -1$; which leads to $S = 0$; which in turn yields

$\left(\text{by Result 5 and Result 6 with } b = \frac{1}{\sqrt{2}}\right)$ the solutions,

$x = 2K\pi \pm \pi/4; \; K \in \mathbb{Z}$.

Lastly, we examine Possibility 2, and we apply Result 10, again with $\kappa = -\sqrt{2}$ and $\lambda = \sqrt{2}$. A necessary and sufficient condition is $f(-\sqrt{2})f(\sqrt{2}) < 0$;

$\Leftrightarrow [4 - \sqrt{2} + (\sqrt{2} - 1)m] \cdot [4 + \sqrt{2} - (\sqrt{2} + 1)m] < 0$; and by factoring out the leading coefficient from each linear factor in m we obtain,

$(\sqrt{2} - 1)(\sqrt{2} + 1) \cdot \left[m + \dfrac{4 - \sqrt{2}}{\sqrt{2} - 1}\right]\left[\dfrac{4 + \sqrt{2}}{1 + \sqrt{2}} - m\right] < 0$; and multiplying by -1 gives,

$= 1$

$\left[m + \dfrac{4 - \sqrt{2}}{\sqrt{2} - 1}\right]\left[m - \left(\dfrac{4 + \sqrt{2}}{1 + \sqrt{2}}\right)\right] > 0 \Leftrightarrow \left\{m > \dfrac{4 + \sqrt{2}}{1 + \sqrt{2}} \text{ or } m < -\left(\dfrac{4 - \sqrt{2}}{\sqrt{2} - 1}\right)\right\}$ (11)

Since m is an integer, (11) is equivalent to
$(m \geq 3) \text{ or } (m \leq -7)$ (12)

And since $m = -n$, (12) is equivalent to
$(n \leq -3) \text{ or } (n \geq 7)$ (13)

The two zeros of the trinomial $f(t) = t^2 + (m-1)t + (2-m)$ are the real numbers

$\dfrac{-(m-1) \pm \sqrt{m^2 + 2m - 7}}{2}$; in terms of n, these two reals are $r_1 = \dfrac{n + 1 + \sqrt{n^2 - 2n - 7}}{2}$

and $r_2 = \dfrac{n + 1 - \sqrt{n^2 - 2n - 7}}{2}$.

We leave it to the reader to verify that when $n \leq -3$, it is $r_1$ which falls inside $(-\sqrt{2}, \sqrt{2})$, while $r_2$ lies outside $[-\sqrt{2}, \sqrt{2}]$. When $n \geq 7$ on the other hand, $r_2$ lies in the open interval $(-\sqrt{2}, \sqrt{2})$, while $r_1$ falls outside $[-\sqrt{2}, \sqrt{2}]$.

Therefore, for $n \leq -3$, $S = r_1$; while for $n \geq 7$, $S = r_2$.

After that we solve for x in accordance with Result 5 and Result 6 (iv). We formalize our findings as follows.

---

**Theorem 1**

Consider the one-variable equation $F(x) = n$, over the set of real numbers $\Re$; where $F(x) = \sin x + \cos x + \sec x + \csc x + \tan x + \cot x$ and n an integer.

1) If $n = -1, 0, 1, 2, 3, 4, 5,$ or $6$; the above equation has no solution in $\Re$.

2) If $n = -2$, all the solutions to the above equation are given by $x = 2K\pi \pm \pi/4$; $K \in Z$.

3) If $n \leq -3$; the real number $r_1 = \dfrac{n + 1 + \sqrt{n^2 - 2n - 7}}{2}$, falls inside the open interval $(-\sqrt{2}, \sqrt{2})$. Let $\varphi_1$ be the unique radian angle measure such that $0 < \varphi_1 < \pi$ and

$\cos\varphi_1 = \dfrac{r_1}{\sqrt{2}}.$ Then all the solutions to the above equation are given by

$x = 2K\pi + \dfrac{\pi}{4} \pm \varphi_1; K \in Z.$

4) If $n \geq 7$, the real number $r_2 = \dfrac{n+1-\sqrt{n^2-2n-7}}{2}$, falls within the open interval $\left(-\sqrt{2}, \sqrt{2}\right)$. Let $\varphi_2$ be the unique radian angle measure such that $0 < \varphi_2 < \pi$ and $\cos\varphi_2 = \dfrac{r_2}{\sqrt{2}}.$ Then all the solutions to the above equation are given by

$x = 2K\pi + \dfrac{\pi}{4} \pm \varphi_2; K \in Z.$

**Remark 1** We invite the reader to verify that $r_1 < 0$, and so, accordingly $\pi/2 < \varphi_1 < \pi$; while $r_2 > 0$, and thus, $0 < \varphi_2 < \pi/2$.